\documentstyle[aps]{revtex}

\begin{document}
\title{Role of Mathematics in Physical Sciences\thanks{%
Based on the talk given at the first national seminar on mathematics
organized by Pakistan Mathematical Society}}
\author{RIAZUDDIN}
\address{National Centre for Physics\\
Quaid-I-Azam University\\
Islamabad 45320, Pakistan}
\date{February 2004}
\maketitle

\begin{abstract}
The role of mathematics in physical sciences is discussed, particularly how
higher mathematics found applications in empirical problems. Several
examples are given to illustrate this role.

(Report No.: NCP-QAU/2004-002)
\end{abstract}

Before I discuss the role of mathematics in physical sciences, let me first
define: What is Mathematics? And What is Science?

\section{What is Mathematics?}

Mathematics is a fusion of skillful operations with concepts and rules
invented just for this purpose \cite{wigner1960}. The principal emphasis is
on invention of concepts, which go beyond those contained in the axioms.
``These are defined with a view of permitting ingenious logical operations
which appeal to our aesthetic sense both as operations and also in their
results of great generality and simplicity'' \cite{wigner1960}. Without
concepts a mathematician would not go far. These may or may not be suggested
by the actual world. ``Mathematics is independent of material objects. In
mathematics the word 'exist' can have only one meaning; it signifies
exemption from contradiction \cite{poincare}.'' In fact most of the advanced
concepts in mathematics are those on which a mathematician can demonstrate
his ingenuity and sense of beauty. Take for example ``complex numbers''; the
introduction of which cannot be suggested by physical observations. A
mathematician's interest in complex numbers lies in that many beautiful
theorems in analytical functions owe their origin to the introduction of
complex numbers. It so happened that much later complex numbers became
essential in the formulation of quantum mechanics where they are not a
calculational trick of applied mathematics. Indeed ``Mathematics can not be
defined without acknowledging its most obvious feature: namely, that it is
interesting'' \cite{polanyi}. It is appropriate to mention Cambridge
Mathematician Godfrey Hardy and his book: ``A Mathematician's Apology'' \cite
{hardy}, the message of which is that pure mathematics is the only kind of
mathematics worthy of respect. He wrote ``A mathematician, like a painter or
a poet is a maker of patterns. If his patterns are more permanent than
theirs, it is because they are made in the ideas. The mathematician's
patterns like the painter's or the poet's must be beautiful; the ideas, like
the colours or the words, must fit together in a harmonious way. Beauty is
the first test, there is no place in the world for ugly mathematics.'' As we
shall see the role which mathematics plays in physical sciences, where one
is concerned to understand basic mysteries of nature, has also beauty of its
own and a source of joy and excitement to its practitioners. Indeed ``a
scientist worthy of name, above all a mathematician, experiences in his work
the same impression as an artist. His pleasure is as great and of the same
nature'' \cite{poincare}.

\section{What is Science (Physics)?}

Science is a fusion of philosophical thinking, which supplies concepts and
skilled crafts, which supply tools. The two are intimately connected.
Concepts are needed to explain old things (in the form of empirical data) in
new ways. Tools are needed to discover new things that have to be explained
(in terms of concepts) or to discover things predicted by a concept-driven
theory so as to verify or discard that theory. The concept driven
revolutions have been rare. Taking quantum mechanics as a prime example of a
concept-driven revolution, Thomas Kuhn in his book ``The Structure of
Scientific revolution'', has listed, in addition to quantum mechanics, only
six major concept driven revolutions in the last 500 years, associated with
the names of Copernicus, Newton, Darwin, Maxwell, Freud and Einstein.
According to F.J. Dyson \cite{dyson}, during the same period, there have
been about twenty tool-driven revolutions, some in physics itself, but
mostly in biology and astronomy, using tools created by physics. Physics has
had great success in creating new tools that have started revolution in
biology, computer science, engineering, astronomy and medicines. Two prime
examples are \cite{dyson} the Galilean revolution resulting from the use of
the telescope in astronomy and the Crick-Watson revolution (1950) resulting
from the use of X-ray crystallography to determine the structure of DNA in
biology. Another tool driven revolution having a great impact on society was
based in the invention of transistor resulting in the advent of computers
and memory banks in the 1960's. Electronic data processing and simulation
revolutionized every branch of science, increasing the power of scientific
theories to interpret and predict new phenomena. Computers, becoming cheaper
and smaller, have become personal and are used for variety of purposes, from
toys to highly sophisticated scientific work. They have revolutionized the
communication, the mode of information and finance.

\section{Role of Mathematics in Physical Theories:}

First there is a mundane role which is to facilitate for the physicists the
numerical calculation of certain constants or the integration of certain
differential equations. Mathematics, does however, play a more sovereign
role in which we will be concerned and bring out how higher mathematics
found applications in subtle empirical problems.

\subsection{The laws of nature are written in the language of mathematics.}

(This is attributed to Galilleo, more than 300 years ago). ``All laws are
deduced from experiment, but to enunciate them, a special language is
needful, ordinary language is too poor, it is besides too vague, to express
relations so delicate, so rich and so precise'' \cite{poincare}. Let us
discuss some examples:

\begin{itemize}
\item[A):]  The basic axioms of Newtonian mathematical physics is stated in
the preface of the first edition of the Principia: rational mechanics ought
to address ``motion'' with the same precision as geometry handles the size
and shape of idealized objects. The association of ``motion'' (particularly
the change in motion) with ``mathematics'' was a stroke of genius. The
mathematical language in which it was formulated contained the concept of
second derivative - not a very immediate concept. The act of writing down a
fundamental law is a rather singular and rare event. It is a miracle that in
spite of the baffling complexity in the world, certain regularities in the
events could be discovered. A monumental example of such a law is Newton's
law of gravitation - a single law which explained everything from planetary
motion to the terrestrial motion of pendulums and which appears simple to
the mathematicians and which proved accurate beyond all reasonable
expectations but still it is a law of limited scope.

\item[B):]  The concepts of modern physics are abstract. ``Mathematics is
the tool specially suited for dealing with abstract concepts of any kind and
there is no limit in its power in this field'' [6]. In this context let us
consider two of the great theories of the last century: Relativity and
Quantum Mechanics, both of which involve mathematics of transformations.
This is because the important quantities in nature appear as the invariant
or having simple transformation properties under these transformations. Let
us consider them one by one:

\begin{enumerate}
\item[i):]  General Theory of Relativity:

Einstein gave a new concept of gravity. Gravity cannot be switched off at
will. Einstein argued that because of its permanency, gravity must be
related to some intrinsic feature of space-time. He identified this feature
as the geometry of space-time $-$ only that that this geometry is unusual.
Existence of matter causes the fabric of space-time to warp somewhat like
the effect of a bowling ball placed on foam. Such distortion to the fabric
of space-time transmits the force of gravity from one place to another.
Gravity resides in the curvature of space-time. The geometry which describes
curved spaces is known as Reimann geometry.

\item[ii):]  Quantum Mechanics:

There are two basic concepts in quantum mechanics: States and Observables.
The states, which have no classical analogue, are vectors in Hilbert space.
The observables are dynamical variables, which although appear in classical
mechanics, are treated in quantum mechanics as hermitian operators on state
vectors. Let us also remind ourselves that Hilbert space of quantum
mechanics is complex with a hermitian scalar product and as such the use of
complex numbers is necessary in the formulation of laws of quantum mechanics.
\end{enumerate}
\end{itemize}

In many cases mathematical concepts were independently developed by the
physicist and recognized then as having been conceived before by the
mathematician. Quantum mechanics is a good example of this where Dirac
invented his own mathematics in his formulation of quantum mechanics.
Einstein, on the other hand, recognised Remiann Geometry as tailor-made for
implementing his view of gravitational force.

\subsection{Mathematical Symmetry and Analogies}

Let us consider this role of mathematics by discussing some examples:

\begin{itemize}
\item[a):]  Maxwell's equations:

The laws of electrodynamics are described by following equations which
expressed all known facts at the time Maxwell began his work: 
\begin{eqnarray*}
{\bf \nabla }\times {\bf H}=\frac{4\pi }c{\bf J} \\
{\bf \nabla }\times {\bf E}+\frac 1c\frac{\partial {\bf B}}{\partial t}=0 \\
{\bf \nabla }\cdot {\bf D}=4\pi \rho \\
{\bf \nabla }\cdot {\bf B}=0
\end{eqnarray*}
In the absence of sources ($\rho =0$, ${\bf J=0}$) Maxwell noticed that the
first two equations lack symmetry as ${\bf \nabla }\times {\bf H}=0 $ while $%
{\bf \nabla }\times {\bf E}=-\left( 1/c\right) \partial {\bf B}/\partial t$.
Maxwell removed this lack of symmetry by modifying the first equation to 
\[
{\bf \nabla }\times {\bf H}=\frac 1c\frac{\partial {\bf D}}{\partial t} 
\]
It was not a new experiment, which came to invalidate the equations. But in
looking at them under a new perspective, Maxwell saw that the equations
become more symmetrical when $\left( 1/c\right) \partial {\bf D}/\partial t$
is added. In that Maxwell was twenty years ahead of experiment since his ``a
priori'' views awaited twenty years for an experimental verification. He
formulated these views because he ``was profoundly steeped in the sense of
mathematical symmetry'' \cite{poincare}. Maxwell unified electricity and
magnetism and as a result the electromagnetic radiation in the form of
light, radiowaves and X-rays provide many of the conveniences of modern life 
$-$ lights, television, telephones etc. Furthermore the requirement of
mutual compatibility of Newtonian Mechanics and Maxwellian Electrodynamics
leads to the foundation of special theory of relativity.

\item[b):]  Dirac Equation

Dirac combined special theory of relativity with quantum mechanics, which
resulted in his famous equation for electron (or any fundamental fermion
having spin 1/2). The Schrodinger equation in non-relativistic quantum
mechanics does not satisfy the basic requirement of relativity, namely that
space and time must be treated on equal footing. This is because it involves
first order time derivative and second other space derivatives. Thus Dirac
looked for an equation which is linear in time and space derivatives. In
order to satisfy this requirement, an analogy with Maxwell's equations,
written earlier, which are Lorentz invariant, may be useful. Maxwell's
equations are first order in time and space derivatives. Moreover, the
vector potential $A_\mu $ to which electric and magnetic fields $E$ and $B$
can be related satisfies in free space the second order wave equation: 
\[
\frac 1{c^2}\frac{\partial ^2}{\partial t^2}A_\mu -\nabla ^2A_\mu =0 
\]
Notice that $E$ and $B$, which satisfy first order equation, have more
components than $A_\mu $ which satisfies second order equation. One should
expect the increase in the number of components as a price one has to pay
for the first order equations. Note also that the above second order
equation is also satisfied by each of the components of $E$ and $B$. The
above considerations suggest that the most general equation we can write is 
\[
\frac 1c\frac{\partial \psi _l}{\partial t}+\sum_{n=1}^N\left( {\bf \alpha }%
\right) _{l\,n}\cdot {\bf \nabla }\psi _n+\frac{imc}\hbar \sum_n\rho
_{l\,n}\psi _n=0 
\]
where ${\bf \alpha }\equiv \left( \alpha _x\text{, }\alpha _y\text{, }\alpha
_z\right) =\alpha _j$, $j=1,2,3$ and $l,n=1,\cdots ,N,$ $N$ being the number
of components we have for the state function $\psi $. Further the
requirement that each component of $\psi $ satisfies a second order equation
(just as $E$ and $B$ do) 
\[
\frac 1{c^2}\frac{\partial ^2\psi }{\partial t^2}-\nabla ^2\psi +\frac{m^2c^2%
}{\hbar ^2}\psi =0 
\]
shows that four quantities $\alpha _x$, $\alpha _y$, $\alpha _z$ and $\beta $
anticommute in pairs and their square is unity. They cannot be numbers and
can be expressed in terms of matrices. One can show that $N$ must be even
and if in the interest of simplicity we require the representation to have
as low a rank as possible, we need to go to $4\times 4$ matrices. In this
way Dirac obtained his equation 
\[
i\hbar \frac{\partial \psi }{\partial t}=H\psi =\left( -c{\bf \alpha }\cdot
c\hbar {\bf \nabla }+\beta mc^2\right) \psi 
\]
where $\psi $ has four components and transforms under Lorentz
transformation in a way different from a vector and is called Dirac spinor.
The Dirac equation has profound consequences: it naturally comes out that
particle it represents has spin 1/2; antimatter must exist, to each particle
there is an antiparticle, and gave a new meaning to vacuum in the
microscopic world. Both Maxwell's equations and Dirac equation give much
more than what was put in, purely from mathematical symmetry and analogies.
Moreover they reveal to a physicist the hidden symmetry of things in making
him see them in a new way.

\item[c):]  Symmetry Principles and Group Theory:

Until the twentieth century group theory played a little role in theoretical
physics. This had a background. ``Most classical physicists expected
infinitesimal analysis to be the natural mathematics for all of physics,
with priority accorded to the differential equations of Newtonian mechanics.
It was a widely held tacit assumption that this must be the way mathematics
enters microphysics \cite{tisza}''. Partly it may also be due to new and
unfamiliar mathematics comprising group theory. As Wigner noted, ``There is
a great reluctance among physicists towards accepting group theoretical
arguments'' \cite{wigner1959}. Wigner's early experience in X-ray
crystallography led him to a programme of applying the theory of group
representations to atomic and molecular spectra as well as nuclear physics.
He also gave the infinite unitary representations of the Lorentz group and
laid foundation both for the application of group theory to quantum
mechanics and for the role of symmetry in microphysics.

In the last three decades Lie groups and Lie algebras played a major role in
applying symmetry principles in containing allowable dynamical laws mainly
for the reason that in many cases in subatomic phenomena, the dynamical laws
were not known a priori.

Let us consider a group element of a continuous special unitary group $SU(N)$
$\exp [i\alpha _aT_a]$ where $\alpha _a$ is a set of real parameters $\left(
a=1,\ldots ,N\right) $, $\det \alpha =1$ and $T_a$ are generators of the
symmetry group satisfying the commutation relations 
\[
\left[ T^a,T^b\right] =if^{abc}T^c 
\]
and the Jacobi identity 
\[
\left[ T^a,\left[ T^b,T^c\right] \right] +\left[ T^b,\left[ T^c,T^a\right]
\right] +\left[ T^c,\left[ T^a,T^b\right] \right] =0 
\]
where fabc are structure constants, antisymmetric in $a,b,c$. The above
equations define the Lie algebra associated with the group SU(N). The group
SU(3) has been used in the classification of physical ``particles'' which
exist in multi-plets and can be regarded as belonging to irreducible
representations of SU(3). The unitary groups have also been used in more
fundamental way when we consider local gauge symmetries where the parameters 
$\alpha _a$ are functions of space-time. Here an object transforms as 
\[
\psi \left( x\right) \rightarrow \exp [i\alpha _a\left( x\right) T_a]\psi
\left( x\right) 
\]
Under an infinitesimal transformation, the derivative 
\[
dx^\mu \partial _\mu \psi \left( x\right) =\left[ \psi \left( x+dx\right)
-\psi \left( x\right) \right] 
\]
does not make sense since $\psi \left( x+dx\right) $ transforms differently
from $\psi \left( x\right) $ under the group. It is necessary to introduce a
vector field $A^\mu $, called gauge field in physics. In the language of
differential geometry $A^\mu $ forms a connection and defines a parallel
displacement of geometrical objects belonging to representation spaces of
the group \cite{itzyksonz}. A parallel displacement of $\psi $ from $x$ to $%
x+dx$ is defined through 
\begin{eqnarray*}
\psi _u\left( x\right) =U\left( x+dx,x,A\right) \psi \left( x\right) \\
=\left[ 1+idx^\mu A_\mu ^aT^a\right] \psi \left( x\right)
\end{eqnarray*}
The definition of covariant derivative then naturally follows 
\begin{eqnarray*}
dx^\mu D_\mu \psi \left( x\right) =\psi \left( x+dx\right) -\psi _u\left(
x\right) \\
=dx^\mu \left( \partial _\mu -iA_\mu ^aT^a\right) \psi \left( x\right)
\end{eqnarray*}
This concept unifies fundamental particles with fundamental forces through
which particles interact; those forces are mediated by the quantas of the
gauge field A. As C.N. Yang stated: ``Symmetry dictates interaction''.
Further, as seen above that the gauge symmetry is based on a sophisticated
geometrical concept, gives it a deep and beautiful foundation.
\end{itemize}

\section{Modern Trends: Quantum Geometry:}

We have two great theories of the last century: the quantum mechanics and
the theory of relativity. The two theories have their roots in mutually
exclusive groups of phenomena. Quantum Mechanics provides a theoretical
framework for understanding the universe on the smallest of scales:
molecules, atoms and all the way to subatomic particles like electrons and
quarks. General relativity provides a theoretical framework for
understanding the universe on the largest of scales: stars, galaxies,
clusters of galaxies, and beyond to the immense expanse of the universe
itself. The two theories operate with different mathematical concepts -
infinite dimensional Hilbert space and the four dimensional Riemann space,
respectively. In most situations their union is not even required. This is
because in most situations as mentioned above the domains (like atoms and
their constituents) in which quantum mechanics is interested and domains
like (status and galasions) in which general theory of relativity is
relevant have no overlap. There are, however, situations where both theories
become relevant. For instance, in a black hole an enormous mass is crumped
to a very small size. At the moment of big bang the whole of the universe
erupted from a microscopic nugget, compared to which even the grain of sand
looks enormous. These are domains that are tiny and yet incredibly massive
and as such require Quantum Theory and General Theory of Relativity
simultaneously. Until recently the two theories could not be united i.e. no
mathematical formulation exists to which both of these theories are
approximations.

It turns out that to achieve this one needs (i) a higher dimensional
space-time for getting both general coordinate invariance of the Einstein
and the Yang-Mills gauge transformations (mentioned earlier) corresponding
to internal degrees of freedom (ii) supersymmetry to avoid tachyons (i.e.
the particles which move faster term the speed of light) and for taming
infinities (iii) going beyond point field theory i.e. the most fundamental
entities are not point-like but extended one dimensional objects. (This too
helps to tame the infinities). The above three ingredients are incorporated
in superstring theory. It naturally contains a massless spin 2 particle
which could be identified by graviton, the mediator of gravitational
interaction just as spin 1 photon is a mediator of electromagnetic
interaction.

First a words about supersymmetry. Supersymmetry incorporates boson-fermion
symmetry. Such theories predict a new kind of matter in the form of
supersymmetric partners of all observed elementary particles. The
observation of these partners would provide the first experimental evidence
for supersymmetry. But there is so far no experimental evidence for such
particles. The experimental situation will become clear in about 5 years
time when the world's largest accelerator being developed at CERN, Geneva
become operational. The mathematics of supersymmetry involves use of
Clifford algebra and Grassman numbers which unlike ordinary numbers,
anticommute. It turned-out that dynamics of superstring theory can be
formulated in 10-dimensional space-time: four familiar space-time dimensions
and six extra dimensions. The extra-spatial dimensions of string theory are
to be ``crumped'' up (to the size of Planck length 10$^{-35}$ m) in a
particular class of 6-dimensional geometrical shapes known as Calabi$-$Yau
shapes. The mathematics of Calabi$-$Yau shapes is studied in a field called
Algebraic Geometry - a relatively new field that combines algebra and
geometry. Towards the end of last century, it led to some of the crowning
achievements of pure mathematics, including the solution of Fermat's last
theorm, Mordell conjecture and the Weil conjectures \cite{griffiths}. It is
now being used by string theorists leading to a new branch of Physics and
Mathematics, which may be called Quantum Geometry \cite{green}.

Let me end this article by quoting from E. P. Wigner's thought provoking
article ``Unreasonable effectiveness of Mathematics in the Natural
Sciences'':

\begin{quote}
``The miracle of the appropriateness of the language of mathematics for the
formulation of the laws of physics is a wonderful gift which we neither
understand nor deserve. We should be grateful for it and hope that it will
remain valid in future research and that it will extend, for better or for
worse, to our pleasure, even though perhaps also to our bafflement, to wide
branches of learning''.
\end{quote}

\underline{{\bf Acknowledgement}}: This work was supported by a grant from
Pakistan Council for Science and Technology.

\end{document}